\documentclass[reqno,11pt]{amsart}

\usepackage{amsmath}
\usepackage{amscd}
\usepackage{amssymb}
\usepackage{enumerate}
\usepackage{amsfonts}
\usepackage{graphicx}
\usepackage[all]{xy}
\usepackage{mathrsfs}

\DeclareFontEncoding{OT2}{}{}

\numberwithin{equation}{section}

\newcommand{\rar}[1]{\stackrel{#1}{\longrightarrow}}

\newcommand{\om}{\omega} \newcommand{\Om}{\Omega}

\newcommand{\bH}{{\mathbb H}}

\newcommand{\bR}{{\mathbb R}}

\newcommand{\cH}{{\mathcal H}}

\newcommand{\fg}{{\mathfrak g}}

\newcommand{\Hom}{\operatorname{Hom}}

\newcommand{\vol}{\operatorname{vol}}

\newcommand{\tens}{\otimes}
\newcommand{\Wedge}{\bigwedge}

\newcommand{\mbr}{\medbreak}

\newcommand{\bra}{\langle}
\newcommand{\ket}{\rangle}
\newcommand{\pairing}{\bra\cdot,\cdot\ket}
\newcommand{\pair}[2]{\bra #1, #2\ket}

\newtheorem{thm}{Theorem}[section]

\newtheorem{prop}[thm]{Proposition}

\theoremstyle{remark}

\newcommand{\hyk}{hyper-K\"ahler\ }

\newcommand{\bun}{\operatorname{Loc}_G^\circ(X)}

\title[Flat $G$-bundles on hyper-K\"ahler manifolds]{The moduli space
of flat $G$-bundles on a compact hyper-K\"ahler manifold}

\author[M.~Abouzaid and M.~Boyarchenko]{Mohammed Abouzaid
\and Mitya
Boyarchenko
}

\thanks{{\em Authors' address:} Department of
Mathematics, University of Chicago, Chicago, IL 60637 \\
{\em E-mail $($M.A.$)$}: {\tt mabouzai@math.uchicago.edu}
\qquad\qquad {\em E-mail $($M.B.$)$}: {\tt mitya@math.uchicago.edu} \\
The research of M.B. was partially supported by NSF grant
DMS-0401164.}

\begin{document}

%

\maketitle




\section{Introduction}\label{s:intro}

\subsection{} Let $(X,g)$ be a compact connected oriented Riemannian manifold,
and let $G$ be a compact connected Lie group with Lie algebra
$\fg$. We study the moduli space of flat $G$-bundles on $X$, which
can be identified with the quotient
\[
\operatorname{Loc}_G(X) = \Hom(\pi_1(X),G)\bigl/ G,
\]
where $\Hom(\pi_1(X),G)$ denotes the space of all group
homomorphisms from the fundamental group $\pi_1(X)$ into $G$, and
$G$ acts on this space by conjugation. Let $\phi:\pi_1(X)\to G$ be
a homomorphism, and let $[\phi]$ denote the corresponding point of
$\operatorname{Loc}_G(X)$. Note that $\phi$ defines a $\pi_1(X)$-module structure
on $\fg$ via the adjoint representation of $G$. Let us denote this
$\pi_1(X)$-module by $\fg_\phi$, and let $E_\phi$ denote the
associated flat vector bundle on $X$. It is well known (see, e.g.,
\cite{hitchin}) that if $[\phi]$ is a smooth point of $\operatorname{Loc}_G(X)$,
then there is a natural identification of the tangent space
\begin{equation}\label{e:tangent}
T_{[\phi]} \operatorname{Loc}_G(X) \cong H^1_{DR}(X,E_\phi).
\end{equation}
The space on the right is the first de Rham cohomology of $X$ with
coefficients in the flat bundle $E_\phi$.

\mbr

Let $B:\fg\times\fg\to\bR$ be a fixed positive definite
$G$-invariant inner product, which always exists since $G$ is
compact. It defines a positive definite metric on $E_\phi$ which
is compatible with the flat connection. Thus we have a notion of a
{\em harmonic} form on $X$ with coefficients in $E_\phi$
(\cite{demailly}, Chapter VI, \S3.3), and Hodge theory gives an
isomorphism
\begin{equation}\label{e:hodge}
\cH^1(X,E_\phi)\cong H^1_{DR}(X,E_\phi),
\end{equation}
where the space on the left is the space of harmonic $1$-forms on
$X$ with coefficients in $E_\phi$. (Cf. \cite{demailly}, Theorem
VI.3.17. That result is stated for flat {\em Hermitian} vector
bundles; however, the proof works just as well in our situation.)
Let
\[
\pairing : \cH^1(X,E_\phi)\times \cH^1(X,E_\phi) \to C^\infty(X)
\]
denote the obvious pairing obtained by combining the metric on
$\Om^1_X$ induced by $g$ and the metric on $E_\phi$ induced by
$B$. It allows us to define an $L^2$ inner product on
$\cH^1(X,E_\phi)$ by
\[
(s,t)\longmapsto \int_X \pair{s}{t}\cdot \vol_g,
\]
where $\vol_g$ denotes the volume form associated to the metric
$g$. The corresponding metric on (the smooth locus of) $\operatorname{Loc}_G(X)$
induced by the identifications \eqref{e:tangent} and
\eqref{e:hodge} is called the {\em $L^2$ metric on the moduli
space of flat $G$-bundles on $X$}.

\subsection{} Our main result is the following
\begin{thm}\label{t:moduli}
If $X$ is \hyk\!\!, then the $L^2$ metric on the smooth locus of
the moduli space of flat $G$-bundles on $X$ is a \hyk metric.
\end{thm}
This result has been known before in the case $\dim X=4$. See,
e.g., \S7 of \cite{hitchin}, where it is proved using the \hyk
reduction technique. If $\dim X>4$, it is not known to us how to
exhibit $\operatorname{Loc}_G(X)$ as the \hyk reduction of another \hyk manifold.
Thus we have to use a different approach, which was suggested by
Karshon's work \cite{karshon}. Using her results we define three
symplectic forms on the smooth locus $\operatorname{Loc}_G^\circ(X)$ of
$\operatorname{Loc}_G(X)$ which are compatible with the $L^2$ metric, such that
the associated complex structures produce the desired \hyk
structure on $\bun$. The details of the proof are presented in
Section \ref{s:proof}.

\subsection{Acknowledgements} We are greatly indebted to Victor
Ginzburg for drawing our attention to the paper \cite{karshon}, and for making helpful comments about our paper. We
would also like to thank Jean-Pierre Demailly and Carlos Simpson
for answering our questions about Hodge theory.


\section{Proof of Theorem \ref{t:moduli}}\label{s:proof}

\subsection{} Let us define a {\em pseudo-\hyk}manifold
to be a manifold $X$ equipped with a pseudo-Riemannian metric $g$
and an action of the algebra of quaternions $\bH$ on the tangent
bundle $TX$ which is parallel with respect to the Levi-Civita
connection of $g$. Then $X$ is \hyk if and only if $g$ is positive
definite. The basic theory of \hyk manifolds does not use the
positive definiteness condition; in particular, if $X$ is a
pseudo-\hyk manifold, then every imaginary unit quaternion
$I\in\bH$ gives rise to an integrable almost complex structure on
$X$ and a real symplectic form $\om_I$ defined by
$\om_I(v,w)=g(Iv,w)$ for all tangent vectors $v,w\in TX$.

\mbr

We begin with the following result, which is undoubtedly well
known to the experts.
\begin{prop}\label{p:hyper-kahler} Let $X$ be a real manifold
equipped with three symplectic forms $\om_I,\om_J,\om_K$.
\begin{enumerate}[(a)]
\item Assume that $\om_{I,J,K}$, when viewed as isomorphisms $TX\to
T^*X$, satisfy the conditions
\begin{equation}\label{e:conds}
\begin{split}
\om_I^{-1}\circ\om_J &= -\om_J^{-1}\circ\om_I, \\
\om_J^{-1}\circ\om_K &= -\om_K^{-1}\circ\om_J, \\
\om_K^{-1}\circ\om_I &= -\om_I^{-1}\circ\om_K.
\end{split}
\end{equation}
Then the isomorphism $g:=\om_I\circ\om_J^{-1}\circ\om_K:TX\to
T^*X$ corresponds to a (possibly indefinite) nondegenerate metric
on $TX$, and $X$ is pseudo-\hyk with respect to this metric.
\item Consequently, if the following condition holds:
\begin{equation}\label{e:posdef}
\begin{split}
& \text{the isomorphism }
g=\om_I\circ\om_J^{-1}\circ\om_K:TX\rar{} T^*X \\
& \text{ defines
a positive definite metric on } X,
\end{split}
\end{equation}
then $X$ is \hyk.
\end{enumerate}
\end{prop}
\begin{proof}
It clearly suffices to prove (a). For any morphism of vector
bundles $A:TX\to T^*X$, we will denote by $A^*:T^*X\to TX$ the
adjoint map. Using the relations $\om_{I,J,K}^*=-\om_{I,J,K}$ and
\eqref{e:conds}, it is easy to check that $g^*=g$, whence $g$
defines a symmetric bilinear form on $TX$, which is obviously
nondegenerate. Next define $I,J,K:TX\to TX$ by
$\om_I(v,w)=g(Iv,w)$, etc. Again, using \eqref{e:conds}, it is
straightforward to check that these automorphisms satisfy the
quaternionic relations
\begin{equation}\label{e:quaternionic}
I^2=J^2=K^2=IJK=-1.
\end{equation}
Finally, it is well known that if a Riemannian manifold is
equipped with three almost complex structures $I,J,K$ satisfying
\eqref{e:quaternionic}, then the manifold is \hyk if and only if
the associated forms $\om_{I,J,K}$ are closed. For a proof, see,
e.g., \cite{nak-book}, Lemma 3.37. This result does not use the
positive definiteness assumption, thus it also applies in our
situation, since $\om_{I,J,K}$ are closed by assumption. This
completes the proof of the proposition.
\end{proof}

\subsection{} We now recall one of the main results of
\cite{karshon}. If $X$ is a compact connected K\"ahler manifold
and $G$ is as in Section 1, then for every
$\phi\in\Hom(\pi_1(X),G)$, we have a well defined map
\[
L : \cH^*(X,E_\phi) \rar{} \cH^{*+2}(X,E_\phi), \quad
\tau\mapsto\om\wedge\tau,
\]
where $\om$ is the K\"ahler form. The hard Lefschetz theorem holds
in this situation and implies in particular that if $\dim_\bR
X=2n$, then the map
\[
L^{n-1} : \cH^1(X,E_\phi) \rar{} \cH^{2n-1}(X,E_\phi)
\]
is an isomorphism. (Both the fact that $L$ is well defined, i.e.,
takes harmonic forms to harmonic forms, and the fact that the last
map is an isomorphism, follow at once from the corresponding
statements for harmonic forms with coefficients in $\bR$ and the
observation that a differential form on $X$ with coefficients in
$E_\phi$ can be written locally as $\tau=\sum\tau_j\tens e_j$,
where $\{e_j\}$ is a flat local orthonormal frame for $E_\phi$.
Then $\tau$ is harmonic if and only if the coefficients $\tau_j$
are harmonic forms on $X$. This follows trivially from the
definitions of the Laplace operators on $\Om^*_X$ and
$\Om^*_X\tens E_\phi$, see \cite{demailly}, Chapter VI, \S3.3.)

\mbr

With this notation, we have
\begin{prop}[\cite{karshon}, Theorem 5]\label{p:karshon}
The bilinear forms $\varpi_\phi$ on $\cH^1(X,E_\phi)\cong
T_{[\phi]}\bun$ defined by
\[
\varpi_\phi(\tau,\eta) = \int_X \tau\wedge L^{n-1}(\eta)
\]
are nondegenerate, skew-symmetric, and together form a {\em
closed} $2$-form on $\bun$.
\end{prop}

\subsection{} Now let $X$ be a compact connected \hyk manifold
with underlying Riemannian metric $g$. We choose three particular
complex structures $I,J,K$ satisfying the quaternionic relations
\eqref{e:quaternionic}. They give rise to the three corresponding
K\"ahler forms $\om_{I,J,K}$. If $\phi\in\Hom(\pi_1(X),G)$, we let
\[
L_I, L_J, L_K : \cH^*(X,E_\phi) \rar{} \cH^{*+2}(X,E_\phi)
\]
denote the corresponding operators on the cohomology of $X$ with
coefficients in $E_\phi$. As in Proposition \ref{p:karshon}, they
induce three nondegenerate forms
\[
\varpi_{I,\phi},\varpi_{J,\phi},\varpi_{K,\phi} \ :\
\cH^1(X,E_\phi) \times \cH^1(X,E_\phi) \rar{} \bR
\]
which define three $2$-forms $\varpi_{I,J,K}$ on $\bun$ that are
closed by Karshon's result. It is now obvious that our Theorem
\ref{t:moduli} follows from the following more precise

\begin{prop}\label{p:linalg}
The three symplectic forms $\varpi_{I,J,K}$ on $\bun$ satisfy the
relations \eqref{e:conds}. Moreover, the isomorphism
$\varpi_I\circ\varpi_J^{-1}\circ\varpi_K : T\bun\to T^*\bun$
corresponds to the $L^2$ metric on $\bun$.
\end{prop}

\subsection{} The proof of Proposition \ref{p:linalg} reduces
immediately to linear algebra, since it only involves statements
about harmonic forms on $X$. Thus, for simplicity, we will change
our setup and notation as follows. Let $V$ be a \hyk vector space,
i.e., a finite dimensional real vector space equipped with a
positive definite inner product $g$ and three endomorphisms
$I,J,K$ satisfying \eqref{e:quaternionic}. The symplectic forms
$\om_{I,J,K}$ and the operators $L_{I,J,K}:\Wedge^*
V^*\to\Wedge^{*+2}V^*$ are defined in the same way as before. It
is clear that in the proof of Proposition \ref{p:linalg}, the
bundle $E_\phi$ can be ignored completely, since all the relations
that we need to check hold trivially ``in the direction of
$E_\phi$.'' Now $\varpi_{I,J,K}$ become symplectic forms on the
dual space $V^*$, and it is straightforward to check that the
relations \eqref{e:conds} in this situation reduce to
\begin{equation}\label{e:conds2}
\bigl(L_I^{n-1}\bigr)^{-1}\circ L_J^{n-1} =
-\bigl(L_J^{n-1}\bigr)^{-1}\circ L_I^{n-1}, \quad \text{etc.,}
\end{equation}
where $\dim_\bR V=2n$. (Of course, by symmetry, it suffices to
check the first equation only.) We will show that, in fact, the
automorphism $\bigl(L_I^{n-1}\bigr)^{-1}\circ L_J^{n-1}$ of $V^*$
coincides with $K^*:V^*\to V^*$, the adjoint of the automorphism
$K$. By symmetry, it will follow that
$\bigl(L_J^{n-1}\bigr)^{-1}\circ L_I^{n-1}=-K^*$, which will imply
\eqref{e:conds2}. Moreover, these equations will also imply that
$\varpi_I\circ\varpi_J^{-1}\circ\varpi_K:V^*\to V^{**}$
corresponds to the dual metric $g^*$ on $V^*$, and the proof of
Proposition \ref{p:linalg} will be complete.

\mbr

Let $\bigstar:\Wedge^* V^*\to \Wedge^{2n-*}V^*$ denote the Hodge
star operator, see \cite{demailly}, Chapter VI, \S3.1. We observe
that $\bigstar^{-1}\circ L_I^{n-1}=n!\cdot I^*$ on $V^*$. Indeed,
let $v,w\in V^*$. Then, by definition,
\[
g^*\bigl( v,\bigstar^{-1} L_I^{n-1}w \bigr) \cdot \vol_g = v\wedge
L_I^{n-1} w = \om_I^{n-1}\wedge v\wedge w =
\om_I(v^\flat,w^\flat)\cdot \om_I^n,
\]
where $v^\flat\in V$ corresponds to $v$ under the isomorphism
$V\to V^*$ defined by $g$. But it is well known that
$\om_I^n=n!\cdot\vol_g$; on the other hand,
\[
\om_I(v^\flat,w^\flat)=g(Iv^\flat,w^\flat)=g^*(v,I^*w).
\]
This shows that $\bigstar^{-1}\circ L_I^{n-1}=n!\cdot I^*$, and
similarly $\bigstar^{-1}\circ L_J^{n-1}=n!\cdot J^*$, as
endomorphisms of $V^*$. Finally, this implies that
\begin{eqnarray*}
\bigl(L_I^{n-1}\bigr)^{-1}\circ L_J^{n-1} &=&
\bigl(\bigstar^{-1}\circ L_I^{n-1}\bigr)^{-1}\circ
\bigstar^{-1}\circ L_J^{n-1} \\
&=& \bigl( n!\cdot I^* \bigr)^{-1}\circ \bigl(n!\cdot J^*\bigr) \\
&=& \bigl(I^*\bigr)^{-1} \circ J^* = (J\circ I^{-1})^* = K^*,
\end{eqnarray*}
which finishes the proof of Proposition \ref{p:linalg}.


\end{document}